\documentclass[a4paper,10pt]{article}
\usepackage{amsmath,amsfonts,bm}
\usepackage{amsthm}
\newcommand{\R}{\textbf{R}}
\newcommand{\C}{\textbf{C}}
\newcommand{\subres}{\textrm{S}}
\newcommand{\imu}{\bm{\mathit{i}}}
\newcommand{\rank}{\textrm{\upshape rank}}
\newcommand{\diag}{\textrm{diag}}

\newtheorem{proposition}{Proposition}
\newtheorem{lemma}{Lemma}

\begin{document}
%%%
\begin{center}
  {\Large\bf 
%%% Title
    GPGCD, an Iterative Method for Calculating Approximate GCD of
    Univariate Polynomials,\\
    with the Complex Coefficients
  }
\mbox{}\\[11pt]
{\large\sc
%%% Author(s)  Let the corresponding author(s) indicated by $*$
  Akira Terui
}
\mbox{}\\[5pt] 
%%% Contact address(es)
Graduate School of Pure and Applied Sciences\\
University of Tsukuba\\
Tsukuba, 305-8571 Japan\\
{\tt terui@math.tsukuba.ac.jp}\\[5pt] 

% %%% The address(es) of the corresponding author(s) 
% \footnotetext[1]{Correspondence to: Postal Address(es) and Tel(fax) number(s).}

\end{center} 

\begin{abstract}
  We present an extension of our GPGCD method, an iterative method for
  calculating approximate greatest common divisor (GCD) of univariate
  polynomials, to polynomials with the complex coefficients.  For a
  given pair of polynomials and a degree, our algorithm finds a pair
  of polynomials which has a GCD of the given degree and whose
  coefficients are perturbed from those in the original inputs, making
  the perturbations as small as possible, along with the GCD.  In our
  GPGCD method, the problem of approximate GCD is transfered to a
  constrained minimization problem, then solved with a so-called
  modified Newton method, which is a generalization of the
  gradient-projection method, by searching the solution iteratively.
  While our original method is designed for polynomials with the real
  coefficients, we extend it to accept polynomials with the complex
  coefficients in this paper.
\end{abstract}

\section{Introduction}

For algebraic computations on polynomials and matrices, approximate
algebraic algorithms are attracting broad range of attentions
recently.  These algorithms take inputs with some ``noise'' such as
polynomials with floating-point number coefficients with rounding
errors, or more practical errors such as measurement errors, then,
with minimal changes on the inputs, seek a meaningful answer that
reflect desired property of the input, such as a common factor of a
given degree.  By this characteristic, approximate algebraic
algorithms are expected to be applicable to more wide range of
problems, especially those to which exact algebraic algorithms were
not applicable.

As an approximate algebraic algorithm, we consider calculating the
approximate greatest common divisor (GCD) of univariate polynomials,
such that, for a given pair of polynomials and a degree $d$, finding a
pair of polynomials which has a GCD of degree $d$ and whose
coefficients are perturbations from those in the original inputs, with
making the perturbations as small as possible, along with the GCD.
This problem has been extensively studied with various approaches
including the Euclidean method on the polynomial remainder sequence
(PRS) (\cite{bec-lab1998b}, \cite{sas-nod89}, \cite{sch1985}), the
singular value decomposition (SVD) of the Sylvester matrix
(\cite{cor-gia-tra-wat1995}, \cite{emi-gal-lom1997}), the QR
factorization of the Sylvester matrix or its displacements
(\cite{cor-wat-zhi2004}, \cite{zar-ma-fai2000}, \cite{zhi2003}),
Pad\'e approximation \cite{pan2001b}, optimization strategies
(\cite{chi-cor-cor1998}, \cite{kal-yan-zhi2006},
\cite{kal-yan-zhi2007}, \cite{kar-lak1998}, \cite{zen0000a}).
Furthermore, stable methods for ill-conditioned problems have been
discussed (\cite{cor-wat-zhi2004}, \cite{ohs-sug-tor1997},
\cite{san-sas2007}).

Among methods in the above, we focus our attention on optimization
strategies.  Already proposed algorithms utilize iterative methods
including the Levenberg-Marquardt method \cite{chi-cor-cor1998}, the
Gauss-Newton method \cite{zen0000a} and the structured total least
norm (STLN) method (\cite{kal-yan-zhi2006}, \cite{kal-yan-zhi2007}).
Among them, STLN-based methods have shown good performance calculating
approximate GCD with sufficiently small perturbations efficiently.

In this paper, we discuss an extension of the GPGCD method, proposed
by the present author \cite{ter2009}, an iterative method with
transferring the original approximate GCD problem into a constrained
optimization problem, then solving it by a so-called modified Newton
method \cite{tan1980}, which is a generalization of the
gradient-projection method \cite{ros1961}.  In the previous paper
\cite{ter2009}, we have shown that our method calculates approximate
GCD with perturbations as small as those calculated by the STLN-based
methods and with significantly better efficiency than theirs.  While
our original method accepts polynomials with the real coefficients as
inputs and outputs in the previous paper, we extend it to handle
polynomials with the complex coefficients in more generalized settings
in this paper.

The rest part of the paper is organized as follows.  In
Section~\ref{sec:formulation}, we transform the approximate GCD
problem into a constrained minimization problem for the case with the
complex coefficients.
% In
% Section~\ref{sec:gp}, we review the framework of the
% gradient-projection method and a modified Newton method.
In Section~\ref{sec:gpgcd}, we show details for calculating the
approximate GCD, with discussing issues in minimizations.  In
Section~\ref{sec:exp}, we demonstrate performance of our algorithm
with experiments.

\section{Formulation of the Approximate GCD Problem}
\label{sec:formulation}

Let $F(x)$ and $G(x)$ be univariate polynomials of degree $m$ and $n$,
respectively, with the complex coefficients
% , given as
% \begin{equation}
%   \label{eq:FG}
%     F(x) = x^m + f_{m-1} x^{m-1} + \cdots + f_0, \quad
%     G(x) = x^n + g_{n-1} x^{n-1} + \cdots + g_0,
% \end{equation}
%with
and $m\ge n>0$.  We permit $F$ and $G$ be relatively prime in general.
For a given integer $d$ satisfying $n\ge d>0$, let us calculate a
deformation of $F(x)$ and $G(x)$ in the form of
%$
\begin{equation*}
%  \label{eq:FtildeGtilde}
    \tilde{F}(x) = F(x) + \varDelta F(x) = H(x)\cdot \bar{F}(x),
    \quad
    \tilde{G}(x) = G(x) + \varDelta G(x) = H(x)\cdot \bar{G}(x),
\end{equation*}
%$
where $\varDelta F(x)$ and $\varDelta G(x)$ are polynomials with the
complex coefficients, whose degrees do not exceed those of $F(x)$ and
$G(x)$, respectively, $H(x)$ is a polynomial of degree $d$, and
$\bar{F}(x)$ and $\bar{G}(x)$ are pairwise relatively prime.  In this
situation, $H(x)$ is an approximate GCD of $F(x)$ and $G(x)$.  For a
given $d$, we try to minimize $\|\varDelta F(x)\|_2^2 + \|\varDelta
G(x)\|_2^2$ the norm of the deformations.

In the case $\tilde{F}(x)$ and $\tilde{G}(x)$ have a GCD of degree
$d$, then the theory of subresultant tells us that the $(d-1)$-th
subresultant of $\tilde{F}$ and $\tilde{G}$ becomes zero, namely we
have
$
%\[
\subres_{d-1}(\tilde{F},\tilde{G})=0,
%\]
$ where $\subres_k(\tilde{F},\tilde{G})$ denotes the $k$-th
subresultant of $\tilde{F}$ and $\tilde{G}$.  Then, the $(d-1)$-th
subresultant matrix $N_{d-1}(F,G)$,
% \begin{equation}
%   \label{eq:nd-1}
%   \begin{split}
%     N_{d-1}(F,G) &=
%     \begin{pmatrix}
%       f_m    &        &        & g_n    &        &  \\
%       \vdots & \ddots &        & \vdots & \ddots &  \\
%       f_0    &        & f_m    & g_0    &        & g_n \\
%              & \ddots & \vdots &        & \ddots & \vdots \\
%              &        & f_0    &        &        & g_0
%     \end{pmatrix},
%     \\
%     &\qquad\;\; \underbrace{\hspace{1.7cm}}_{n-d+1}
%     \;\;\; \underbrace{\hspace{1.5cm}}_{m-d+1}
%   \end{split}
% \end{equation}
where the $k$-th subresultant matrix $N_k(F,G)$ is a submatrix of the
Sylvester matrix $N(F,G)$ by taking the left $n-k$ columns of
coefficients of $F$ and the left $m-k$ columns of coefficients of $G$,
has a kernel of dimension equal to $1$.  Thus, there exist
polynomials $A(x),B(x)\in\C[x]$ satisfying
\begin{equation}
  \label{eq:abcond}
  A\tilde{F}+B\tilde{G}=0,
\end{equation}
with $\deg(A)<n-d$ and $\deg(B)<m-d$ and $A(x)$ and $B(x)$ are
relatively prime. Therefore, for the given $F(x)$, $G(x)$ and $d$, our
problem is to find $\varDelta F(x)$, $\varDelta G(x)$, $A(x)$ and
$B(x)$ satisfying Eq.\ (\ref{eq:abcond}) with making $\|\varDelta
F\|_2^2+\|\varDelta G\|_2^2$ as small as possible.

Let us assume that $F(x)$ and $G(x)$ are represented as
\begin{equation*}
%  \label{eq:FG-complex}
  \begin{split}
    F(x) &= (f_{m,1}+f_{m,2}\imu) x^m + \cdots + (f_{0,1}+f_{0,2}\imu)
    = F_{\text{Re}}(x) + \imu F_{\text{Im}}(x),
    \\
    G(x) &= (g_{n,1}+g_{n,2}\imu) x^n + \cdots + (g_{0,1}+g_{0,2}\imu)
    = G_{\text{Re}}(x) + \imu G_{\text{Im}}(x),
  \end{split}
\end{equation*}
where $f_{j,1}$, $g_{j,1}$, $f_{j,2}$, $g_{j,2}$ are the real numbers
and $\imu$ is the imaginary unit, and $F_{\text{Re}}(x)$ and
$G_{\text{Re}}(x)$ represent the real part of $F(x)$ and $G(x)$,
respectively, while $F_{\text{Im}}(x)$ and $G_{\text{Im}}(x)$
represent the imaginary part of $F(x)$ and $G(x)$, respectively.
Furthermore, we represent $\tilde{F}(x)$, $\tilde{G}(x)$, $A(x)$ and
$B(x)$ with the complex coefficients as
\begin{equation}
  \label{eq:fgab-complex}
  \begin{split}
    \tilde{F}(x) &= (\tilde{f}_{m,1}+\tilde{f}_{m,2}\imu) x^m +\cdots+
    (\tilde{f}_{0,1}+\tilde{f}_{0,2}\imu) x^0
    = \tilde{F}_{\text{Re}}(x) + \imu \tilde{F}_{\text{Im}}(x),
    \\
    \tilde{G}(x) &= (\tilde{g}_{n,1}+\tilde{g}_{n,2}\imu) x^n +\cdots+
    (\tilde{g}_0x^0+\tilde{g}_{0,2}\imu) x^0
    = \tilde{G}_{\text{Re}}(x) + \imu \tilde{G}_{\text{Im}}(x),
    \\
    A(x) &= (a_{n-d,1}+a_{n-d,2}\imu) x^{n-d} +\cdots+
    (a_{0,1}+a_{0,2}\imu)x^0
    = A_{\text{Re}}(x) + \imu A_{\text{Im}}(x),
    \\ 
    B(x) &= (b_{m-d,1}+b_{m-d,2}\imu) x^{m-d} +\cdots+
    (b_{0,1}+b_{0,2}\imu)x^0
    = B_{\text{Re}}(x) + \imu B_{\text{Im}}(x), 
  \end{split}
\end{equation}
respectively, where $\tilde{f}_{j,1}$, $\tilde{f}_{j,2}$,
$\tilde{g}_{j,1}$, $\tilde{g}_{j,2}$, $a_{j,1}$, $a_{j,2}$, $b_{j,1}$,
$b_{j,2}$ are the real numbers, and, as in above,
$\tilde{F}_{\text{Re}}(x)$, $\tilde{G}_{\text{Re}}(x)$,
$A_{\text{Re}}(x)$ and $B_{\text{Re}}(x)$ represent the real part of
$\tilde{F}(x)$, $\tilde{G}(x)$, $A(x)$ and $B(x)$, respectively, while
$\tilde{F}_{\text{Im}}(x)$, $\tilde{G}_{\text{Im}}(x)$,
$A_{\text{Im}}(x)$ and $B_{\text{Im}}(x)$ represent the imaginary part
of $\tilde{F}(x)$, $\tilde{G}(x)$, $A(x)$ and $B(x)$, respectively.

For the objective function, $\|\varDelta F\|_2^2+\|\varDelta G\|_2^2$
becomes as
\begin{equation}
  \label{eq:objective-complex}
  \sum_{j=0}^m[(\tilde{f}_{j,1}-f_{j,1})^2 +
  (\tilde{f}_{j,2}-f_{j,2})^2] +
  \sum_{j=0}^n[(\tilde{g}_{j,1}-g_{j,1})^2 +
  (\tilde{g}_{j,2}-g_{j,2})^2].
\end{equation}

For the constraint, Eq.\ (\ref{eq:abcond}) becomes as
\begin{multline}
  \label{eq:abcond2-complex}
%   \begin{pmatrix}
%     \tilde{f}_{m,1}+\tilde{f}_{m,2}\imu & & &
%     \tilde{g}_{n,1}+\tilde{g}_{n,2}\imu & &  \\
%     \vdots      & \ddots &             & \vdots      & \ddots &  \\
%     \tilde{f}_{0,1}+\tilde{f}_{0,2}\imu & &
%     \tilde{f}_{m,1}+\tilde{f}_{m,2}\imu &
%     \tilde{g}_{0,1}+\tilde{g}_{0,2}\imu & &
%     \tilde{g}_{n,1}+\tilde{g}_{n,2}\imu \\
%     & \ddots & \vdots &    & \ddots      & \vdots \\
%     & & \tilde{f}_{0,1}+\tilde{f}_{0,2}\imu & & &
%     \tilde{g}_{0,1}+\tilde{g}_{0,2}\imu 
%   \end{pmatrix}
%   \\
%   \times
  N_{d-1}(\tilde{F},\tilde{G})\cdot
  {}^t
  \begin{pmatrix}
    a_{n-d,1}+a_{n-d,2}\imu, \ldots, a_{0,1}+a_{0,2}\imu,
    b_{m-d,1}+b_{m-d,2}\imu, \ldots, b_{0,1}+b_{0,2}\imu
  \end{pmatrix}
  \\
  =
  \bm{0}
  .
\end{multline}
By expressing the subresultant matrix and the column vector in
\eqref{eq:abcond2-complex} separated into the real and the complex parts,
respectively, we express \eqref{eq:abcond2-complex} as
\begin{gather}
%\begin{equation}
  \label{eq:abcond2-complex-2}
  (N_1+N_2\imu)(\bm{v}_1+\bm{v}_2\imu)=\bm{0},
  \\
%\end{equation}
%with
%\begin{equation}
  \label{eq:n1n2v1v2}
  \begin{split}
    N_1 &= N_{d-1}(\tilde{F}_{\text{Re}}(x), \tilde{G}_{\text{Re}}(x)),
%     N_1 &=
%     \begin{pmatrix}
%       \tilde{f}_{m,1} & & & \tilde{g}_{n,1} & &  \\
%       \vdots      & \ddots &             & \vdots      & \ddots &  \\
%       \tilde{f}_{0,1} & & \tilde{f}_{m,1} & \tilde{g}_{0,1} & &
%       \tilde{g}_{n,1} \\ 
%       & \ddots & \vdots &    & \ddots      & \vdots \\
%       & & \tilde{f}_{0,1} & & & \tilde{g}_{0,1}
%     \end{pmatrix},
%     \\
    \quad
    N_2 = N_{d-1}(\tilde{F}_{\text{Im}}(x), \tilde{G}_{\text{Im}}(x)),
%     N_2 &=
%     \begin{pmatrix}
%       \tilde{f}_{m,2} & & & \tilde{g}_{n,2} & &  \\
%       \vdots      & \ddots &             & \vdots      & \ddots &  \\
%       \tilde{f}_{0,2} & & \tilde{f}_{m,2} & \tilde{g}_{0,2} & &
%       \tilde{g}_{n,2} \\ 
%       & \ddots & \vdots &    & \ddots      & \vdots \\
%       & & \tilde{f}_{0,2} & & & \tilde{g}_{0,2}
%     \end{pmatrix},
    \\
    \bm{v}_1 &= {}^t(a_{n-d,1},\ldots,a_{0,1},b_{m-d,1},\ldots,b_{0,1}),
%    \\
    \quad
    \bm{v}_2 
%    &
    = {}^t(a_{n-d,2},\ldots,a_{0,2},b_{m-d,2},\ldots,b_{0,2}).
  \end{split}
%\end{equation}
\end{gather}
We can expand the left-hand-side of Eq.\ \eqref{eq:abcond2-complex-2} as
$
%\begin{equation*}
  (N_1+N_2\imu)(\bm{v}_1+\bm{v}_2\imu)=
  (N_1\bm{v}_1-N_2\bm{v}_2)+\imu(N_1\bm{v}_2+N_2\bm{v}_1),
%\end{equation*}
$
thus, Eq.\ \eqref{eq:abcond2-complex-2} is equivalent to a system of
equations:
$
%\begin{equation*}
  N_1\bm{v}_1-N_2\bm{v}_2=\bm{0},
%  \quad
  N_1\bm{v}_2+N_2\bm{v}_1=\bm{0},
%\end{equation*}
$
which is expressed as
\begin{equation}
  \label{eq:abcond2-complex-3}
  \begin{pmatrix}
    N_1 & -N_2 \\
    N_2 & N_1
  \end{pmatrix}
  \begin{pmatrix}
    \bm{v}_1 \\ \bm{v}_2
  \end{pmatrix}
  =
  \bm{0}.
\end{equation}

Furthermore, we add another constraint for the coefficient of $A(x)$
and $B(x)$ as
\begin{multline}
  \label{eq:abconstraint-complex}
  \|A(x)\|_2^2+\|B(x)\|_2^2
  =(a_{n-d,1}^2+\cdots+a_{0,1}^2)+(b_{m-d,1}^2+\cdots+b_{0,1}^2)
  \\
  +(a_{n-d,2}^2+\cdots+a_{0,2}^2)+(b_{m-d,2}^2+\cdots+b_{0,2}^2)-1=0,
\end{multline}
which can be expressed together with \eqref{eq:abcond2-complex-3} as
\begin{equation}
  \label{eq:abcond2-complex-4}
  \begin{pmatrix}
    {}^t\bm{v}_1 & {}^t\bm{v}_2 & -1 \\
    N_1 & -N_2 & \bm{0} \\
    N_2 & N_1 & \bm{0}
  \end{pmatrix}
  \begin{pmatrix}
    \bm{v}_1 \\ \bm{v}_2 \\ 1
  \end{pmatrix}
  =
  \bm{0},
\end{equation}
where Eq.\ \eqref{eq:abconstraint-complex} has been put on the top of
Eq.\ \eqref{eq:abcond2-complex-3}.  Note that, in Eq.\
\eqref{eq:abcond2-complex-4}, we have total of $2(m+n-d+1)+1$
equations in the coefficients of polynomials in
\eqref{eq:fgab-complex} as a constraint, with the $j$-th row of which
is expressed as $q_j=0$.

Now, we substitute the variables
\begin{multline}
  \label{eq:var0-complex}
  (
  \tilde{f}_{m,1},\ldots,\tilde{f}_{0,1},
  \tilde{g}_{n,1},\ldots,\tilde{g}_{0,1},
  \tilde{f}_{m,2},\ldots,\tilde{f}_{0,2},
  \tilde{g}_{n,2},\ldots,\tilde{g}_{0,2},
  \\
  a_{n-d,1},\ldots,a_{0,1},
  b_{m-d,1},\ldots,b_{0,1},
  a_{n-d,2},\ldots,a_{0,2},
  b_{m-d,2},\ldots,b_{0,2}
  ),
\end{multline}
as $\bm{x}=(x_1,\ldots,x_{4(m+n-d+2)})$, then Eq.\
\eqref{eq:objective-complex} and \eqref{eq:abcond2-complex-4} become
as
\begin{multline}
  \label{eq:objective-complex2}
  f(\bm{x})=
  (x_1-f_{m,1})^2+\cdots+(x_{m+1}-f_{0,1})^2
  +
  (x_{m+2}-g_{n,1})^2+\cdots+(x_{m+n+2}-g_{0,1})^2
  \\
  +
  (x_{m+n+3}-f_{m,2})^2+\cdots+(x_{2m+n+3}-f_{0,2})^2
  \\
  +
  (x_{2m+n+4}-g_{n,2})^2+\cdots+(x_{2(m+n+2)}-g_{0,2})^2,
\end{multline}
\begin{equation}
  \label{eq:constraint-complex2}
  \bm{q}(\bm{x})=
  {}^t(q_1(\bm{x}), \ldots, q_{2(m+n-d+1)+1}(\bm{x}))
  =
  \bm{0},
\end{equation}
respectively.  Therefore, the problem of finding an approximate GCD
can be formulated as a constrained minimization problem of finding a
minimizer of the objective function $f(\bm{x})$ in Eq.\
(\ref{eq:objective-complex2}), subject to $\bm{q}(\bm{x})=\bm{0}$ in Eq.\
(\ref{eq:constraint-complex2}).

\section{The Algorithm for Approximate GCD}
\label{sec:gpgcd}

We calculate an approximate GCD by solving the constrained
minimization problem \eqref{eq:objective-complex2},
\eqref{eq:constraint-complex2} with the gradient projection method by
Rosen \cite{ros1961} (whose initials become the name of our GPGCD
method) or a modified Newton method by Tanabe \cite{tan1980} (for
review, see the author's previous paper \cite{ter2009}).  Our
preceding experiments \cite[Section 5.1]{ter2009} have shown that a
modified Newton method was more efficient than the original gradient
projection method while the both methods have shown almost the same
convergence property, thus we adopt a modified Newton method in this
paper.

In applying a modified Newton method to the approximate GCD problem,
we discuss issues in the construction of the algorithm in detail, such
as
\begin{itemize}
\item Representation of the Jacobian matrix $J_{\bm{g}}(\bm{x})$ and
  certifying that $J_{\bm{g}}(\bm{x})$ has full rank
  (Section \ref{sec:jacobianmat}),
\item Setting the initial values (Section~\ref{sec:init}),
\item Regarding the minimization problem as the minimum distance
  problem (Section~\ref{sec:mindistance}),
\item Calculating the actual GCD and correcting the coefficients of
  $\tilde{F}$ and $\tilde{G}$ (Section~\ref{sec:correct}),
\end{itemize}
as follows.
% In this paper, we discuss the complex coefficient case;
% for the real coefficient case, see the author's previous paper
% \cite{ter2009}.
% After presenting the algorithm, we give a modification for
% preserving monicity, and end this section with examples.

\subsection{Representation and the rank of the  Jacobian Matrix}
\label{sec:jacobianmat}

For a polynomial $P(x)\in\C[x]$ represented as
%\[
$
P(x) = p_nx^n+\cdots+p_0x^0,
$
%\]
let $C_k(P)$ be a complex $(n+k,k+1)$ matrix defined as
\[
\begin{array}{ccl}
  C_k(P) & = & 
  \begin{pmatrix}
    p_n & & \\
    \vdots & \ddots & \\
    p_0 & & p_n \\
    & \ddots & \vdots \\
    & & p_0
  \end{pmatrix}
  .
  \\[-2mm]
  & & 
  \hspace{3mm}
  \underbrace{\hspace{17mm}}_{k+1}
\end{array}
\]
For co-factors $A(x)$ and $B(x)$ in \eqref{eq:fgab-complex},
% consider
% matrices $C_m(A)$ and $C_n(B)$ and express them as the sum of matrices
% consisting of the real and the imaginary parts of whose elements,
% respectively, as
% \begin{equation}
%   \begin{split}
%     C_m(A) &=
%     \begin{pmatrix}
%       a_{n-d,1} & & \\
%       \vdots & \ddots & \\
%       a_{0,1} & & a_{n-d,1} \\
%       & \ddots & \vdots \\
%       & & a_{0,1}
%     \end{pmatrix}
%     +\imu
%     \begin{pmatrix}
%       a_{n-d,2} & & \\
%       \vdots & \ddots & \\
%       a_{0,2} & & a_{n-d,2} \\
%       & \ddots & \vdots \\
%       & & a_{0,2}
%     \end{pmatrix}
%     \\
%     &= C_m(A)_1 + \imu C_m(A)_2,
%     \\
%     C_n(B) &=
%     \begin{pmatrix}
%       b_{m-d,1} & & \\
%       \vdots & \ddots & \\
%       b_{0,1} & & b_{m-d,1} \\
%       & \ddots & \vdots \\
%       & & b_{0,1}
%     \end{pmatrix}
%     +\imu
%     \begin{pmatrix}
%       b_{m-d,2} & & \\
%       \vdots & \ddots & \\
%       b_{0,2} & & b_{m-d,2} \\
%       & \ddots & \vdots \\
%       & & b_{0,2}
%     \end{pmatrix}
%     \\
%     &= C_n(B)_1 + \imu C_n(B)_2,
%   \end{split}
% \end{equation}
% respectively, and define
define matrices $A_1$ and $A_2$ as
\begin{equation}
  \label{eq:a1a2}
  A_1 = [C_m(A_{\text{Re}}(x))\; C_n(B_{\text{Re}}(x))],
  \quad
  A_2 = [C_m(A_{\text{Im}}(x))\; C_n(B_{\text{Im}}(x))].
%   \begin{split}
%     A_1 &= [C_m(A)_1\; C_n(B)_1]
%     =
%     \begin{pmatrix}
%       a_{n-d,1} & & & b_{m-d,1} & & \\
%       \vdots & \ddots & & \vdots & \ddots & \\
%       a_{0,1} & & a_{n-d,1} & b_{0,1} & & b_{m-d,1} \\
%       & \ddots & \vdots & & \ddots & \vdots \\
%       & & a_{0,1} & & & b_{0,1}
%     \end{pmatrix}
%     ,
%     \\
%     A_2 &= [C_m(A)_2\; C_n(B)_2]
%     =
%     \begin{pmatrix}
%       a_{n-d,2} & & & b_{m-d,2} & & \\
%       \vdots & \ddots & & \vdots & \ddots & \\
%       a_{0,2} & & a_{n-d,2} & b_{0,2} & & b_{m-d,2} \\
%       & \ddots & \vdots & & \ddots & \vdots \\
%       & & a_{0,2} & & & b_{0,2}
%     \end{pmatrix}
%     .
%   \end{split}
\end{equation}
(Note that $A_1$ and $A_2$ are matrices of the real numbers of
$m+n-d+1$ rows and $m+n+2$ columns.)  Then, by the definition of the
constraint \eqref{eq:constraint-complex2}, we have the Jacobian matrix
$J_{\bm{g}}(\bm{x})$ (with the original notation of variables
\eqref{eq:var0-complex} for $\bm{x}$) as
\begin{equation*}
%  \label{eq:jacobian-complex}
  J_{\bm{g}}(\bm{x}) =
  \begin{pmatrix}
    \bm{0} & \bm{0} & 2\cdot{}^t\bm{v}_1 & 2\cdot{}^t\bm{v}_2 \\
    A_1 & -A_2 & N_1 & -N_2 \\
    A_2 & A_1 & N_2 & N_1
  \end{pmatrix}
  ,
\end{equation*}
with $A_1$ and $A_2$ as in \eqref{eq:a1a2} and $N_1$, $N_2$,
$\bm{v}_1$ and $\bm{v}_2$ as in \eqref{eq:n1n2v1v2}, respectively,
which can be easily constructed in every iteration.
%in Algorithms~\ref{alg:gp} and \ref{alg:gpnewton}.

In executing
%Algorithm~\ref{alg:gp} or \ref{alg:gpnewton},
iterations, we need to keep that $J_{\bm{g}}(\bm{x})$ has full rank:
otherwise,
% we cannot correctly calculate $(J_{\bm{g}}(\bm{x}))^+$ (in
% Algorithm~\ref{alg:gp}) or the matrix in \eqref{eq:modnewton} becomes
% singular (in Algorithm~\ref{alg:gpnewton}) thus
we are unable to
decide proper search direction.  For this requirement, we have the
following observations.
\begin{proposition}
  \label{prop:fullrank}
  Let $\bm{x}^*\in V_{\bm{g}}$ be any feasible point satisfying Eq.\
  \eqref{eq:constraint-complex2}.  Then, if the corresponding
  polynomials do not have a GCD whose degree exceeds $d$, then
  $J_{\bm{g}}(\bm{x}^*)$ has full rank.
\end{proposition}
\begin{proof}
  Let $\bm{x}^*=
  (
  \tilde{f}_{m,1},\ldots,\tilde{f}_{0,1},
  \tilde{g}_{n,1},\ldots,\tilde{g}_{0,1},
  \tilde{f}_{m,2},\ldots,\tilde{f}_{0,2},
  \tilde{g}_{n,2},\ldots,\tilde{g}_{0,2},
  a_{n-d,1},\ldots,a_{0,1},
  $
  $
  b_{m-d,1},\ldots,b_{0,1},
  a_{n-d,2},\ldots,a_{0,2},
  b_{m-d,2},\ldots,b_{0,2}
  )
  $ with its polynomial representation \linebreak
  expressed as in
  (\ref{eq:fgab-complex}) (note that this assumption permits the
  polynomials $\tilde{F}(x)$ and $\tilde{G}(x)$ to be relatively prime
  in general).  To verify our claim, we show that we have
  $\rank(J_{\bm{g}}(\bm{x}^*))=2(m+n-d+1)+1$.
  % as in (\ref{eq:rankj}).
  Let us express $J_{\bm{g}}(\bm{x}^*)=
  \begin{pmatrix}
    J_\mathrm{L} \mid J_\mathrm{R}
  \end{pmatrix}
  $, where $J_\mathrm{L}$ and $J_\mathrm{R}$ are column blocks
  expressed as
  \[
  J_\mathrm{L} =
  \begin{pmatrix}
    \bm{0} & \bm{0} \\
    A_1 & -A_2 \\
    A_2 & A_1 
  \end{pmatrix}
  ,
  \quad
  J_\mathrm{R} =
  \begin{pmatrix}
    2\cdot{}^t\bm{v}_1 & 2\cdot{}^t\bm{v}_2 \\
    N_1 & -N_2 \\
    N_2 & N_1
  \end{pmatrix}
  ,
  \]
  respectively.  Then, we have the following lemma.
  \begin{lemma}
    \label{lem:fullrank}
    We have $\rank(J_\mathrm{L})=2(m+n-d+1)$.
  \end{lemma}
  \begin{proof}
    For $A_1=[C_m(A)_1\; C_n(B)_1]$, let $\overline{C_m(A)_1}$ be the
    right $m-d$ columns of $C_m(A)_1$ and $\overline{C_n(B)_1}$ be the
    right $n-d$ columns of $C_n(B)_1$.  Then, we see that the bottom
    $m+n-2d$ rows of the matrix $\bar{C}=[\overline{C_m(A)_1}\;
    \overline{C_n(B)_1}]$ is equal to the matrix consisting of the
    real part of the elements of $N(A,B)$, the Sylvester matrix of
    $A(x)$ and $B(x)$.  By the assumption, polynomials $A(x)$ and
    $B(x)$ are relatively prime, and there exist no nonzero elements
    in $\bar{C}$ except for the bottom $m+n-2d$ rows, thus we have
    $\rank(\bar{C})=m+n-2d$.

    By the structure of $\bar{C}$ and the lower triangular structure
    of $C_m(A)_1$ and $C_n(B)_1$, we can take the left $d+1$ columns
    of $C_m(A)_1$ or $C_n(B)_1$ satisfying linear independence along
    with $\bar{C}$, which implies that there exist a nonsingular
    square matrix $T$ of order $m+n+2$ satisfying
    \begin{equation}
      \label{eq:tr1}
      A_1 T=R,
    \end{equation}
    where $R$ is a lower triangular matrix, thus we have
    $\rank(A_1)=\rank(R)=m+n-d+1$.

    Furthermore, by using $T$ and $R$ in \eqref{eq:tr1}, we have
    \begin{equation}
      \label{eq:tr2}
      \begin{pmatrix}
        \bm{0} & \bm{0} \\
        A_1 & -A_2 \\
        A_2 & A_1
      \end{pmatrix}
      \begin{pmatrix}
        T & \bm{0} \\
        \bm{0} & T
      \end{pmatrix}
      =
      \begin{pmatrix}
        \bm{0} & \bm{0} \\
        R & -A_2T \\
        A_2T & R
      \end{pmatrix}
      ,
    \end{equation}
    followed by a suitable transformation on columns on the matrix in
    the right-hand-side of \eqref{eq:tr2}, we can make $A_2T$ to zero
    matrix, which implies that
    \[
    \rank(J_\mathrm{L})=\rank
    \left(
      \begin{pmatrix}
        \bm{0} & \bm{0} \\
        R & -A_2T \\
        A_2T & R
      \end{pmatrix}
    \right)
    =2\cdot\rank(R)=2(m+n-d+1).
    \]
    This proves the lemma.
  \end{proof}

  \textsc{Proof of Proposition~\ref{prop:fullrank} (continued).}  By
  the assumptions, we have at least one nonzero coordinate in the top
  row in $J_\mathrm{R}$, while we have no nonzero coordinate in the
  top row in $J_\mathrm{L}$, thus we have
  $\rank(J_{\bm{g}}(\bm{x}))=2(m+n-d+1)+1$, which proves the
  proposition.
\end{proof}

Proposition~\ref{prop:fullrank} says that, so long as the search
direction in the minimization problem satisfies that corresponding
polynomials have a GCD of degree not exceeding $d$, then
$J_{\bm{g}}(\bm{x})$ has full rank, thus we can safely calculate the
next search direction for approximate GCD.

\subsection{Setting the Initial Values}
\label{sec:init}

At the beginning of iterations, we give the initial value $\bm{x}_0$
by using the singular value decomposition (SVD) \cite{dem1997} of
% In
% the case of the real coefficients, we calculate the SVD of the
% $(d-1)$-th subresultant matrix $N_{d-1}(F,G):
% \R^{m+n-2d}\rightarrow\R^{m+n-d}$ in (\ref{eq:nd-1}) (see the author's
% previous paper \cite{ter2009}).  In the complex case, we calculate the
% SVD of 
$
N=
\begin{pmatrix}
  N_1 & -N_2 \\
  N_2 & N_1
\end{pmatrix}
$
in \eqref{eq:abcond2-complex-3} as
$
%\begin{equation}
%  \label{eq:svd-nd-1}
%  \begin{array}{c}
    N =  U\,\Sigma\,{}^tV,
%    \\
    U = (\bm{u}_1,\ldots,\bm{u}_{2(m+n-2d+2)}),
%    \quad
$
$
    \Sigma = \diag(\sigma_1,\ldots,\sigma_{2(m+n-2d+2)}),
%    \\
    V = (\bm{v}_1,\ldots,\bm{v}_{2(m+n-2d+2)}),
%  \end{array}
%\end{equation}
$
with $\bm{u}_j\in\R^{2(m+n-d+1)}$, $\bm{v}_j\in\R^{2(m+n-2d+2)}$, and
$\Sigma=\diag(\sigma_1,\ldots,$ $\sigma_{2(m+n-2d+2)})$ denotes the
diagonal matrix whose the $j$-th diagonal element is $\sigma_j$.  Note
that $U$ and $V$ are orthogonal matrices.  Then, by a property of the
SVD \cite[Theorem~3.3]{dem1997}, the smallest singular value
$\sigma_{2(m+n-2d+2)}$ gives the minimum distance of the image of the
unit sphere $\textrm{S}^{2(m+n-2d+2)-1}$, given as
%\[
$
\textrm{S}^{2(m+n-2d+2)-1}=
\{
\bm{x}\in\R^{2(m+n-2d+2)} \mid \|\bm{x}\|_2=1
\}
,
$
%\]
by $N$, represented as
\linebreak
%\[
$
N\cdot\textrm{S}^{2(m+n-2d+1)-1}=
\{
N\bm{x}\mid \bm{x}\in\R^{2(m+n-2d+2)}, \|\bm{x}\|_2=1
\},
$
%\]
from the origin, along with $\sigma_{2(m+n-2d+2)}\bm{u}_{2(m+n-2d+2)}$
as its coordinates.
%By (\ref{eq:svd-nd-1}),
Thus, we have
%\[
$
N\cdot \bm{v}_{2(m+n-2d+2)} = \sigma_{2(m+n-2d+2)}\bm{u}_{2(m+n-2d+2)}.
$
%\]
% thus  $\bm{v}_{2(m+n-2d+2)}$ represents the coefficients
% of $A(x)$ and $B(x)$: let
For
$\bm{v}_{m+n-2d}={}^t(\bar{a}_{n-d},\ldots,\bar{a}_0,\bar{b}_{n-d},\ldots,\bar{b}_0)$,
let
\linebreak
$\bar{A}(x) = \bar{a}_{n-d}x^{n-d}+\cdots+\bar{a}_0x^0$ and
$\bar{B}(x) = \bar{b}_{m-d}x^{m-d}+\cdots+\bar{b}_0x^0$.
% \[
% \begin{split}
%   \bm{v}_{2(m+n-2d+2)} &=
%   {}^t(
%   \bar{a}_{n-d,1},\ldots,\bar{a}_{0,1},
%   \bar{b}_{n-d,1},\ldots,\bar{b}_{0,1},
%   \bar{a}_{n-d,2},\ldots,\bar{a}_{0,2},
%   \bar{b}_{n-d,2},\ldots,\bar{b}_{0,2}
%   ),
%   \\
%   \bar{A}(x) &= (\bar{a}_{n-d,1}+\bar{a}_{n-d,2}\imu)x^{n-d}
%   +\cdots+
%   (\bar{a}_{0,1}+\bar{a}_{0,2}\imu) x^0, \\
%   \bar{B}(x) &= (\bar{b}_{m-d,1}+\bar{b}_{m-d,2}\imu)x^{m-d}
%   +\cdots+
%   (\bar{b}_{0,1}+\bar{b}_{0,2}\imu) x^0.
% \end{split}
% \]
Then, $\bar{A}(x)$ and $\bar{B}(x)$ give the least norm of $AF+BG$
satisfying $\|A\|_2^2+\|B\|_2^2=1$ by putting $A(x)=\bar{A}(x)$ and
$B(x)=\bar{B}(x)$ in \eqref{eq:fgab-complex}. 

Therefore, we admit the coefficients of $F$, $G$, $\bar{A}$ and
$\bar{B}$ as the initial values of the iterations as
\begin{multline*}
%  \label{eq:appgcd-init}
  \bm{x}_0 = 
  (
  f_{m,1},\ldots,f_{0,1},
  g_{n,1},\ldots,g_{0,1},
  f_{m,2},\ldots,f_{0,2},
  g_{n,2},\ldots,g_{0,2},
  \\
  \bar{a}_{n-d,1},\ldots,\bar{a}_{0,1},
  \bar{b}_{n-d,1},\ldots,\bar{b}_{0,1},
  \bar{a}_{n-d,2},\ldots,\bar{a}_{0,2},
  \bar{b}_{n-d,2},\ldots,\bar{b}_{0,2}
  ).
\end{multline*}

\subsection{Regarding the Minimization Problem as the Minimum Distance
  (Least Squares) Problem} 
\label{sec:mindistance}

Since we have the object function $f$ as in
(\ref{eq:objective-complex2}), we have
\begin{multline*}
  \nabla f(\bm{x})=
  2\cdot {}^t(
  x_1-f_{m,1},\ldots,x_{m+1}-f_{0,1},
  x_{m+2}-g_{n,1},\ldots,x_{m+n+2}-g_{0,1},
  \\
  x_{m+n+3}-f_{m,2},\ldots,x_{2m+n+3}-f_{0,2},
  x_{2m+n+4}-g_{n,2},\ldots,x_{2(m+n+2)}-g_{0,2},
  0,\ldots,0).
\end{multline*}
 However, we can regard our problem as finding a point
 $\bm{x}\in V_{\bm{g}}$ which has the minimum distance to the initial
 point $\bm{x}_0$ with respect to the
 $(x_1,\ldots,x_{2(m+n+2)})$-coordinates which correspond to the
 coefficients in $F(x)$ and $G(x)$.
%  Therefore, in the gradient projection method at $\bm{x}\in
%  V_{\bm{g}}$, the projection of $-\nabla f(\bm{x})$ in
%  (\ref{eq:projection}) should be the projection of
% \begin{multline}
%   {}^t(
%   x_1-f_{m,1},\ldots,x_{m+1}-f_{0,1},
%   x_{m+2}-g_{n,1},\ldots,x_{m+n+2}-g_{0,1},
%   \\
%   x_{m+n+3}-f_{m,2},\ldots,x_{2m+n+3}-f_{0,2},
%   x_{2m+n+4}-g_{n,2},\ldots,x_{2(m+n+2)}-g_{0,2},
%   0,\ldots,0),
% \end{multline}
% onto $T_{\bm{x}}$. This change is equivalent to changing the
% objective function as $\bar{f}(\bm{x})=\frac{1}{2}f(\bm{x})$ then solving
% the minimization problem of $\bar{f}(\bm{x})$, subject to
% $\bm{q}(\bm{x})=\bm{0}$.
 Therefore, as in the real case (see the authors previous paper
 \cite{ter2009}), we change the objective function as
 $\bar{f}(\bm{x})=\frac{1}{2}f(\bm{x})$, then solve the minimization
 problem of $\bar{f}(\bm{x})$, subject to $\bm{q}(\bm{x})=\bm{0}$.

\subsection{Calculating the Actual GCD and Correcting the Deformed
  Polynomials}
\label{sec:correct}

After successful end of the iterations,
% in Algorithms~\ref{alg:gp} or \ref{alg:gpnewton},
we obtain the coefficients of $\tilde{F}(x)$,
$\tilde{G}(x)$, $A(x)$ and $B(x)$ satisfying (\ref{eq:abcond}) with
$A(x)$ and $B(x)$ are relatively prime.  Then, we need to compute the
actual GCD $H(x)$ of $\tilde{F}(x)$ and $\tilde{G}(x)$.  Although $H$
can be calculated as the quotient of $\tilde{F}$ divided by $B$ or
$\tilde{G}$ divided by $A$, naive polynomial division may cause
numerical errors in the coefficient.  Thus, we calculate the
coefficients of $H$ by the so-called least squares division
\cite{zen0000a}, followed by correcting the coefficients in
$\tilde{F}$ and $\tilde{G}$ by using the calculated $H$, as follows.

For polynomials $\tilde{F}$, $\tilde{G}$, $A$ and $B$
represented as in \eqref{eq:fgab-complex} and $H$ represented as
%\[
$
H(x)=(h_{d,1}+h_{d,2}\imu)x^d+\cdots+(h_{0,1}+h_{0,2}\imu)x^0,
$
%\]
solve the equations $HB=\tilde{F}$ and $HA=\tilde{G}$ with respect to
$H$ as solving the least squares problems of linear systems
\begin{align}
  \label{eq:hsystem1}
  C_d(A)\,{}^t(h_{d,1}+h_{d,2}\imu,\ldots,h_{0,1}+h_{0,2}\imu) &=
  {}^t(\tilde{g}_{n,1}+\tilde{g}_{n,2}\imu,\ldots,
  \tilde{g}_{0,1}+\tilde{g}_{0,2}\imu), \\
  \label{eq:hsystem2}
  C_d(B)\,{}^t(h_{d,1}+h_{d,2}\imu,\ldots,h_{0,1}+h_{0,2}\imu) &=
  {}^t(\tilde{f}_{m,1}+\tilde{f}_{m,2}\imu,\ldots,
  \tilde{f}_{0,1}+\tilde{f}_{0,2}\imu),
\end{align}
respectively.  Then, we transfer the linear systems
\eqref{eq:hsystem1} and \eqref{eq:hsystem2}, as follows.  For
\eqref{eq:hsystem2}, let us express the matrices and vectors as the
sum of the real and the imaginary part of which, respectively, as
$
%\begin{gather*}
  C_d(B) = B_1+\imu B_2,
%  \\
  {}^t(h_{d,1}+h_{d,2}\imu,\ldots,h_{0,1}+h_{0,2}\imu) =
  \bm{h}_1+\imu\bm{h}_2,
%  \\
  {}^t(\tilde{f}_{m,1}+\tilde{f}_{m,2}\imu,\ldots,
  \tilde{f}_{0,1}+\tilde{f}_{0,2}\imu) =
  \bm{f}_1+\imu\bm{f}_2.
%\end{gather*}
$
Then, \eqref{eq:hsystem2} is expressed as
\begin{equation}
  \label{eq:hsystem2-2}
  (B_1+\imu B_2)(\bm{h}_1+\imu\bm{h}_2)=(\bm{f}_1+\imu\bm{f}_2).
\end{equation}
By equating the real and the imaginary parts in Eq.\
\eqref{eq:hsystem2-2}, respectively, we have
%\[
$
(B_1\bm{h}_1-B_2\bm{h}_2)=\bm{f}_1,
%\quad
(B_1\bm{h}_2+B_2\bm{h}_1)=\bm{f}_2,
$
%\]
or
\begin{equation}
  \label{eq:hsystem2-3}
  \begin{pmatrix}
    B_1 & -B_2 \\
    B_2 & B_1
  \end{pmatrix}
  \begin{pmatrix}
    \bm{h}_1 \\ \bm{h}_2
  \end{pmatrix}
  =
  \begin{pmatrix}
    \bm{f}_1 \\ \bm{f}_2
  \end{pmatrix}
  .
\end{equation}
Thus, we can calculate the coefficients of $H(x)$ by solving the real
least squares problem \eqref{eq:hsystem2-3}.  We can solve
\eqref{eq:hsystem1} similarly.

Let $H_1(x),H_2(x)\in\C[x]$ be the candidates for the
GCD whose coefficients are calculated as the least squares solutions
of (\ref{eq:hsystem1}) and (\ref{eq:hsystem2}), respectively.  Then,
for $i=1,2$, calculate the norms of the residues as
%\[
$
r_i = \|\tilde{F}-H_iB\|_2^2+\|\tilde{G}-H_iA\|_2^2,
$
%\]
respectively, and set the GCD $H(x)$ be $H_i(x)$ giving the minimum
value of $r_i$.

Finally, for the chosen $H(x)$, correct the coefficients of
$\tilde{F}(x)$ and $\tilde{G}(x)$ as
% \[
% \tilde{F}(x)=H(x)\cdot B(x), \quad
% \tilde{G}(x)=H(x)\cdot A(x),
% \]
$\tilde{F}(x)=H(x)\cdot B(x)$, $\tilde{G}(x)=H(x)\cdot A(x)$,
respectively.

\section{Experiments}
\label{sec:exp}

We have implemented the GPGCD algorithm for polynomials with the
complex coefficients
% (Algorithm~\ref{alg:gpgcd})
on the computer algebra system Maple and compared its performance with
a method based on the structured total least norm (STLN) method
\cite{kal-yan-zhi2006}
% and carried out the following tests:
% \begin{itemize}
% \item Comparison of performance of the gradient-projection method
%   (Algorithm~\ref{alg:gp}) and a modified Newton method
%   (Algorithm~\ref{alg:gpnewton}),
% \item Comparison of performance of the GPGCD method with a method based on
%   the structured total least norm (STLN) method
%   \cite{kal-yan-zhi2006}, 
% \end{itemize}
for randomly generated polynomials with approximate GCD.
The tests have been carried out on Intel Core2 Duo Mobile Processor
T7400 (in Apple MacBook ``Mid-2007'' model) at $2.16$ GHz with RAM
2GB, under MacOS X 10.5.

In the tests, we have generated random polynomials with GCD then added
noise, as follows.  First, we have generated a pair of monic
polynomials $F_0(x)$ and $G_0(x)$ of degrees $m$ and $n$,
respectively, with the GCD of degree $d$.  The GCD and the prime parts
of degrees $m-d$ and $n-d$ are generated as monic polynomials and with
random coefficients $c\in[-10,10]$ of floating-point numbers.  For
noise, we have generated a pair of polynomials $F_{\mathrm{N}}(x)$ and
$G_{\mathrm{N}}(x)$ of degrees $m-1$ and $n-1$, respectively, with
random coefficients as the same as for $F_0(x)$ and $G_0(x)$.  Then,
we have defined a pair of test polynomials $F(x)$ and $G(x)$ as
%$
\[
%\begin{split}
  F(x) 
%  &
  = F_0(x)+\frac{e_F}{\|F_{\mathrm{N}}(x)\|_2}F_{\mathrm{N}}(x),
%  \\
  \quad
  G(x)
%  &
  = G_0(x)+\frac{e_G}{\|G_{\mathrm{N}}(x)\|_2}G_{\mathrm{N}}(x),
%\end{split}
\]
%$
respectively, scaling the noise such that the $2$-norm of the noise
for $F$ and $G$ is equal to $e_F$ and $e_G$, respectively.  In the
present test, we set $e_F=e_G=0.1$.

% \subsection{Tests on Large Sets of Randomly-generated Polynomials}
% \label{sec:test-appgcd}

In this test, we have compared our implementation against a method
based on the structured total least norm (STLN) method
\cite{kal-yan-zhi2006}, using their implementation (see
Acknowledgments).  In their STLN-based method, we have used the
procedure \verb|C_con_mulpoly| which calculates the approximate GCD of
several polynomials in $\C[x]$.  The tests have been carried out on
Maple 12 with \verb|Digits=15| executing hardware floating-point
arithmetic.  For every example, we have generated $100$ random test
polynomials as in the above.  In executing
%Algorithm~\ref{alg:gpgcd},
a modified Newton method, 
we set
%$u=200$ and
a threshold of the $2$-norm of the ``update'' vector in each
iteration $\varepsilon=1.0\times 10^{-8}$; in
\verb|C_con_mulpoly|, we set the tolerance $e=1.0\times 10^{-8}$.

Table~\ref{tab:appgcd} shows the results of the test: $m$ and $n$
denotes the degree of a pair $F$ and $G$, respectively, and $d$
denotes the degree of approximate GCD.  The columns with ``STLN'' are
the data for the STLN-based method, while those with ``GPGCD'' are the
data for the GPGCD method.
% ``\#Fail'' is the
% number of ``failed'' cases such as: in the STLN-based method, the
% number of iterations exceeds $50$ times (which is the built-in
% threshold in the program), while, in the GPGCD method, the
% perturbation \eqref{eq:perturbation}
% \begin{equation}
%   \label{eq:perturbation}
%   \|\tilde{F}-F\|_2^2+\|\tilde{G}-G\|_2^2
% \end{equation}
% exceeds $1$ (note that, in the GPGCD method, all the iterations have
% converged within far less than $200$ times).  All the other data are
% the average over results for the ``not failed'' cases:
% ``Error'',
% ``\#Iterations'' and ``Time'' are the same as those in
% Table~\ref{tab:gp}, respectively.
``Error'' is the perturbation 
$
%\begin{equation}
%  \label{eq:perturbation}
  \|\tilde{F}-F\|_2^2+\|\tilde{G}-G\|_2^2,
%\end{equation}
$
where ``$ae\!-\!b$'' denotes $a\times 10^{-b}$; ``\#Iterations'' is
the number of iterations; ``Time'' is computing time in seconds. 

We see that, in the most of tests, both methods calculate approximate
GCD with almost the same amount of perturbations, while the GPGCD
method runs much faster than STLN-based method by approximately from
$10$ to $30$ times.  Note that, in contrast to the real coefficient
case \cite{ter2009}, both methods have converged in all the test cases
with the number of iterations and sufficiently small amount of
perturbations as approximately equal to those shown as in
Table~\ref{tab:appgcd}.

\begin{table*}
  \centering
  \begin{tabular}{|c|c|c|c|c|c|c|c|c|}
    \hline
    Ex. & $m,n$ & $d$ & 
    \multicolumn{2}{c|}{Error}  & \multicolumn{2}{c|}{\#Iterations} &
    \multicolumn{2}{c|}{Time (sec.)} \\
    \cline{4-9}
    & & & STLN & GPGCD & STLN & GPGCD & STLN & GPGCD \\
    \hline
    1 & $10,10$ & $5$ & $3.72e\!-\!3$ & $3.72e\!-\!3$ &
    $4.48$ & $4.43$ & $1.79$ & $0.15$\\ 
    \hline
    2 & $20,20$ & $10$ & $4.16e\!-\!3$ & $4.16e\!-\!3$&
    $4.24$ & $4.22$ & $5.88$ & $0.30$ \\ 
    \hline
    3 & $30,30$ & $15$ & $4.33e\!-\!3$ & $4.33e\!-\!3$&
    $4.54$ & $4.48$ & $14.29$ & $0.58$ \\ 
    \hline
    4 & $40,40$ & 20 & $4.48e\!-\!3$ & $4.48e\!-\!3$ &
    $4.08$ & $4.08$ & $24.10$ & $0.88$ \\  
    \hline
    5 & $50,50$ & 25 & $4.63e\!-\!3$ & $4.64e\!-\!3$ &
    $4.05$ & $4.12$ & $39.19$ & $1.36$ \\ 
    \hline
    6 & $60,60$ & 30 & $4.61e\!-\!3$ & $4.61e\!-\!3$ &
    $4.02$ & $4.06$ & $60.48$ & $1.96$ \\ 
    \hline
    7 & $70,70$ & 35 & $4.82e\!-\!3$ & $4.82e\!-\!3$ &
    $3.90$ & $4.02$ & $84.51$ & $2.66$ \\ 
    \hline
    8 & $80,80$ & 40 & $4.84e\!-\!3$ & $4.84e\!-\!3$ &
    $3.88$ & $4.04$ & $116.03$ & $3.65$ \\ 
    \hline
    9 & $90,90$ & 45 & $4.79e\!-\!3$ & $4.79e\!-\!3$ &
    $3.85$ & $4.01$ & $151.27$ & $4.66$ \\ 
    \hline
    10 & $100,100$ & 50 & $4.77e\!-\!3$ & $4.78e\!-\!3$ &
    $3.83$ & $4.06$ & $199.48$ & $6.00$ \\ 
    \hline
  \end{tabular}
  \caption{Test results for large sets of polynomials with approximate
    GCD. See Section~\ref{sec:exp} for details.}
  \label{tab:appgcd}
\end{table*}

\section{Concluding Remarks}

Based on our previous research \cite{ter2009}, we have extended our
GPGCD method for polynomials with the complex coefficients.

Our experiments have shown that, as in the real coefficients case
\cite{ter2009}, our algorithm calculates approximate GCD with
perturbations as small as those calculated by methods based on the
structured total least norm (STLN) method, while our method has shown
significantly better performance over the STLN-based methods in its
speed, by approximately up to $30$ times, which seems to be
sufficiently practical for inputs of low or moderate degrees.  This
result shows that, in contrast to their \textit{structure preserving}
method, our simple method can achieve accurate and efficient
computation as or more than theirs in calculating approximate GCDs.

Our future research includes theoretical investigation of convergence
properties, investigation for efficient computation in solving a linear
system in each iteration by analysis of the structure of matrices,
generalization of our method to several input polynomials, and so on.

\section*{Acknowledgments}

% If you would like to thank anyone or your grant(s), place your comments here
% and remove the percent signs.

We thank Professor Erich Kaltofen for making their implementation for
computing approximate GCD available on the Internet and providing
experimental results.  We also thank the anonymous reviewers for their
helpful comments.

This research was supported in part by the Ministry of Education,
Culture, Sports, Science and Technology of Japan, under Grant-in-Aid
for Scientific Research (KAKENHI) 19700004.

% BibTeX users please use 
% \bibliographystyle{plain}    
% \bibliography{terui-e}

%\begin{thebibliography}{109}
%Author(s), Article title, Journal, Volume, page numbers, Year.
%Author(s), Book title, page numbers, Publisher, place, Year.
% etc
%\end{thebibliography}

\end{document}